\documentclass[12pt]{article}
\usepackage{fullpage}

\usepackage{xcolor}
\usepackage[normalem]{ulem}

\usepackage{graphicx} 

\usepackage{amsmath,amssymb,amsmath}

\usepackage[square,numbers]{natbib}
\usepackage[final,colorlinks]{hyperref}

\usepackage{url}

\definecolor{darkblue}{cmyk}{1,0,0,0.5}
\definecolor{darkred}{cmyk}{0,1,0,0.5}

\hypersetup{anchorcolor=black,
  citecolor=darkred, filecolor=darkblue,
  menucolor=darkblue,urlcolor=darkblue,linkcolor=darkblue}

\newcommand{\R}{\mathbb{R}}

\newcommand{\Tr}{\mathrm{Tr}}
\newcommand{\tran}{\mathsf{T}}

\graphicspath{{./figures/}}

\title{Local interaction of two systems with saddle-node bifurcations: mutualistic and mixed cases}
\author{Peter Ashwin,\thanks{
Department of Mathematics and Statistics,
University of Exeter,
Exeter EX4 4QF, UK}
\and Claire Postlethwaite\thanks{Department of Mathematics, University of Auckland, Auckland, 1142, New Zealand}
\and Jan Sieber${}^*$}

\date{\today}

\begin{document}

\maketitle

\begin{abstract}
The saddle-node bifurcation is the simplest example of a generic bifurcation in smooth ordinary differential equations, and is associated with the creation or destruction of a pair of equilibria. In this paper we examine the unfolding of the dynamics that occur when two generically coupled systems have simultaneous saddle-node bifurcations. We note that four parameters are required to generically unfold the interactions, and the dynamics are surprisingly complicated relative to the simplicity of a single saddle-node bifurcation. In the unfolding, in addition to saddle-node, Hopf and codimension-two local bifurcations, we also find a variety of global bifurcations, including homoclinic, SNIC, SNICeroclinic and non-central SNIC bifurcations. The latter two are codimension-two bifurcations that occur at the termination of a curve of SNIC bifurcations. A further contribution of this work is the development of numerical continuation techniques for the tracking of these codimension-two bifurcations through parameter space.
\end{abstract}


\newpage

\section{Introduction}

The saddle-node, fold or limit point bifurcation is the simplest example of a generic bifurcation that appears in a smooth nonlinear ordinary differential equation (ODE) with one parameter. The topological normal form on the centre manifold  $x\in \R$ with parameter  
$\lambda\in\R$ is given by
\begin{equation}
\frac{d}{dt}x=\dot{x}=x^2+\lambda
\end{equation}
on varying $\lambda$ through zero. This local bifurcation can be associated with a ``SNIC'' bifurcation where there is creation of periodic orbits in a global bifurcation, for example if $x$ is on a periodic domain. This SNIC bifurcation is present, for example, in
\begin{equation}
\frac{d}{dt}x=\dot{x}=(\cos x+1)+\lambda,
\end{equation}
if we consider $x$ as an angle on the unit circle ($\R\!\!\mod\!(2\pi)$) and vary $\lambda$ through zero.

In this paper, we examine the dynamics that appear on generically unfolding two coupled systems, each with a saddle-node bifurcation at some point in parameter space. Although many bifurcations more complex than this have been analysed in the past, the interaction of two weakly coupled systems, each undergoing a saddle-node bifurcation, motivates the problem studied here, namely bifurcation of an equilibrium whose Jacobian has two zero eigenvalues and geometric multiplicity two.
Checking the linear constraints at an equilibrium undergoing such a bifurcation, it is necessary for the Jacobian to have zero eigenvalues with algebraic multiplicity two. Generically, such a Jacobian will have a nontrivial Jordan block, meaning we have geometric multiplicity one (i.e.~only one eigenvector); this is the condition necessary for a Takens--Bogdanov bifurcation \cite{kuznetsov2004elements}, generic for two parameters. At the interaction of two saddle nodes, we have geometric multiplicity two as well - this implies four parameters are needed to generically unfold an interaction of two saddle nodes. 

Although it is codimension four, it is still of interest in that it appears naturally at the uncoupling limit \cite{perez2019uncoupling} of two saddle-node bifurcations. The interaction of two saddle-node bifurcations was previously analysed in the 1970s as the ``hilltop bifurcation'' \cite{thompson1973general}. 
The generic unfolding of equilibria in the neighbourhood of such a bifurcation with one distinguished parameter and three unfolding parameters is given in \cite{golubitsky1979theory} and \cite[Chapter IX]{golubitsky1985}. However, we are not aware of a systematic attempt to unfold such a bifurcation in a way that respects not only equilibria but also the dynamics near this local bifurcation point. In the context of this bifurcation with global reinjection, interacting saddle-node bifurcations on invariant circles were analysed by Baesens and MacKay \cite{baesens2013interaction} in the ``mutualistic'' case where the linear interactions have the same signs. They highlight a second, more difficult, ``mixed'' case where the linear interactions have opposite signs. We analyse the local bifurcations for the mixed case in detail and so extend their analysis. A further example is provided in \cite{augustsson2024coevolutionary}, who consider the global dynamics of a parametrised system consisting of two coupled saddle-node bifurcations on a torus. They consider the mixed case, but in a symmetric case and do not provide a full local unfolding.

The paper is organised as follows. In Section~\ref{sec:analysis}, we discuss the normal form, exhibit an invariance under rescaling and some parameter/phase space symmetries that help simplify the unfolding. We identify the locus of codimension one Hopf/neutral saddles and saddle-node bifurcations analytically. For the mutualistic case (Section~\ref{sec:mutual}), only saddle-node bifurcations are present. For the mixed case (Section~\ref{sec:mixed}), we also find analytic expressions for codimension-two Takens--Bogdanov, cusp and Bautin bifurcations. The remaining bifurcations need to be found by numerical continuation. We first find the locus of homoclinic and SNIC bifurcations by fixing two parameters ($\alpha$ and $\beta$) and varying two more parameters ($\gamma$ and $\mu$). The curve of SNIC bifurcations terminates at one end in a non-central SNIC bifurcation, and at the other in a heteroclinic involving a saddle node and a hyperbolic equilibrium, also called a SNICeroclinic bifurcation in \cite{nechyporenko2024}. 
We do a complete local bifurcation analysis of the mutualistic and mixed cases using the continuation package COCO~\cite{dankowicz2013recipes}. This includes the development of new numerical continuation schemes (detailed in Appendix~\ref{sec:defining:systems}) to continue the codimension-two SNICeroclinic and non-central SNIC bifurcations. To the best of our knowledge, this is the first time these particular codimension-two bifurcations have been continued numerically.

Section~\ref{sec:transient} turns to implications for transient behaviour. 
In particular we  analyse the phase plane to understand when there are transient trapping regions as well as trapping regions associated with local attractors. 
Finally, Section~\ref{sec:discuss} highlights various issues, including generalizing these local results to understand the bifurcations of mixed case flows on the torus, parallel to the work in \cite{baesens2013interaction} in the mutualistic case.

\section{Normal form analysis}
\label{sec:analysis}

The local (topological) normal form considered in \cite[Eq 3]{baesens2013interaction} is written as follows:
\begin{equation}
     \begin{split}
		\dot{x}&=x^2-\lambda- \gamma +2\alpha y +\epsilon_1 y^2+HOT_1 \\
			\dot{y}&=y^2-\lambda+ \gamma+2\beta x +\epsilon_2x^2+HOT_2 
	\end{split} \label{eq:baesens-normform}
 \end{equation}
with parameters $\alpha$, $\beta$, $\lambda$, $\gamma$, $\epsilon_{1,2}$ and higher order terms indicated by $HOT_{1,2}$. Baesens et al \cite{baesens2013interaction} consider cases where $\alpha$ and $\beta$ are of similar order to each other; more precisely that $\alpha\beta\neq 0$. Motivated by coupled SNIC (saddle-node on invariant circle) bifurcations, they consider this local normal form on a torus, where escape to infinity corresponds to a transient with re-injection. We are also motivated by this, though we do not consider the global dynamics in the paper: the local cases are already quite complicated.

Baesens et al \cite{baesens2013interaction} identify three cases and analyse the mutualistic case in some detail:
\begin{itemize}
    \item The case $\alpha\beta>0$ is called {\em mutualistic} ($\alpha>0$, $\beta>0$ called {\em mutualistic excitatory}
    and $\alpha<0$, $\beta<0$ called {\em mutualistic inhibitory}).
    \item The case $\alpha\beta<0$ is called {\em mixed}.
    \item The case $\alpha\beta=0$ is called {\em degenerate}.
\end{itemize} 
The mixed case is not analysed in \cite{baesens2013interaction} - one of the main results in this paper is to include this case. Note that for differentiable equivalence, the cubic terms in the higher order terms in (\ref{eq:baesens-normform}) cannot be removed \cite{glendinning2022normal}, though we shall follow \cite{baesens2013interaction} in ignoring these terms.

Away from $\alpha\beta=0$, the analysis in~\cite{baesens2013interaction} suggests that one can truncate the $\epsilon$ terms to give the main system that we will consider, namely \cite[Eq 4]{baesens2013interaction}: 
\begin{equation}
     \begin{split}
		\dot{x}&=x^2-\lambda- \gamma+2\alpha y, \\
			\dot{y}&=y^2-\lambda+ \gamma+2\beta x 
	\end{split} \label{eq:odes}
 \end{equation}
where $\lambda$, $\alpha$, $\beta$ and $\gamma$ are parameters. In this section, we give some analytical results on the bifurcation structure of this planar system. We present results discussed in~\cite{baesens2013interaction} for the mutualistic case $\alpha\beta>0$, but also consider the mixed case
where $\alpha\beta<0$.

We show in Section~\ref{sec:hopf} that if we define
\begin{equation}
    \lambda_{{\Tr}0}(\gamma,\alpha,\beta):=\frac{\gamma(\gamma+\alpha^2-\beta^2)}{(\alpha+\beta)^2}
    \label{eq:lambdaTr0}
\end{equation}
then (\ref{eq:odes}) has an equilibrium with trace zero precisely when $\lambda=\lambda_{{\Tr}0}$. It is convenient to write (\ref{eq:odes}) in the form
\begin{equation}
     \begin{split}
		\dot{x}&=x^2-\mu-\lambda_{{\Tr}0}- \gamma+2\alpha y, \\
			\dot{y}&=y^2-\mu-\lambda_{{\Tr}0}+ \gamma+2\beta x 
	\end{split} \label{eq:odesmu}
\end{equation}
where $\mu$, $\alpha$, $\beta$ and $\gamma$ are parameters that we  use in the latter part of this paper. 


\subsection{Rescaling and parameter/phase space symmetries}
\label{sec:scaling}

The system (\ref{eq:odes}) is invariant under the parameter rescaling
\begin{equation}
    (\alpha,\beta,\gamma,\lambda)=(\rho\cos\theta,\rho\sin\theta,\rho^2\hat{\gamma},\rho^2\hat{\lambda})
    \label{eq:params}
\end{equation}
and the time and coordinate rescaling
\[
(x,y,t)=(\rho \hat{x},\rho \hat{y}, \rho^{-1} \hat{t}).
\]
This means that by an isotropic rescaling of phase space and time we can remove one of the parameters. In particular, we can assume without loss of generality that $\alpha^2+\beta^2$ is constant and non-zero and then by scaling of $\rho$ we will recover all possibilities up to rescaling.

Note that the equations remain unchanged under each of the following two transformations of coordinates and variables:
\[
\alpha\rightarrow-\alpha,\quad \beta\rightarrow -\beta,\quad x\rightarrow -x,\quad y\rightarrow -y, \quad t\rightarrow -t
\]
\[
\alpha\rightarrow\beta,\quad \beta\rightarrow \alpha,\quad x\rightarrow y,\quad y\rightarrow x, \quad \gamma\rightarrow -\gamma
\]
In the case $\alpha\beta>0$ we can thus further assume that $\alpha>\beta>0$. When $\alpha\beta<0$, we can assume that $\alpha>0>\beta$ and $\alpha>|\beta|$.
Note that because of the change of direction of time in the first transformation, although the computation of bifurcation curves remains unchanged under this transformation, stabilities of equilibria, and hence criticality of bifurcations and direction of flow in phase portraits will change. 

\subsection{Degeneracies of the normal form}

One degeneracy is in the case $\alpha\beta=0$; although neither $\alpha=0$ nor $\beta=0$ are special for the scaling (\ref{eq:params}), as noted by \cite{baesens2013interaction} they do give rise to a degeneracy of the truncated normal form (\ref{eq:odes}). For $\alpha=0$ this corresponds to the $y$ dynamics being a skew product over the $x$ dynamics ($x$ drives $y$) while for $\beta=0$ the $x$ dynamics is a skew product over the $y$ dynamics ($y$ drives $x$). Complete absence of coupling in one direction is a degeneracy (codimension infinity), such that one needs to take into account the higher order terms such as those with coefficients $\epsilon_{i}$ in (\ref{eq:baesens-normform}), as highlighted in~\cite{baesens2013interaction}.

Another degeneracy appears in the case $\alpha=-\beta$ and $\gamma=0$, namely in this case the normal form (\ref{eq:odes}) becomes 
\begin{equation}
     \begin{split}
		\dot{x}&=x^2-\lambda+2\alpha y, \\
		\dot{y}&=y^2-\lambda-2\alpha x 
	\end{split} \label{eq:reversing}
 \end{equation}
which has a reversing symmetry 
$$
(x,y,t)\mapsto(-y,-x,-t)
$$
with reversing symmetry line $x=-y$. This more subtle degeneracy implies that equilibria in the line of symmetry have zero trace (saddles will be resonant), attractors/repellors must lie off this plane in symmetric pairs, and any periodic orbits crossing the symmetry line will come in families and be neutrally stable. Higher order terms would have to be taken into account to resolve this degeneracy.

For $\alpha=\beta$ and $\gamma=0$ the normal form is
\begin{equation}
     \begin{split}
		\dot{x}&=x^2-\lambda+2\alpha y, \\
		\dot{y}&=y^2-\lambda+2\alpha x 
	\end{split} \label{eq:xysymm}
\end{equation}
which has a symmetry 
$$
(x,y)\mapsto (y,x)
$$
that has an invariant line $x=y$. This symmetry will, in a similar way, be broken by higher order terms. If we consider these bifurcations as associated with vector fields that pass through a neighbourhood of the origin in a similar way to (\ref{eq:odes}) note that addition just of cubic order terms (as in \cite{glendinning2022normal}) may not be enough: we suspect that quartic terms may need to be included to break the degeneracies of this two-dimensional normal form. In this paper, we assume we are far enough away from these degenerate points that the additional terms are not needed: we leave a full exploration of these points to future work.


\subsection{Hopf bifurcations and neutral saddles}
\label{sec:hopf}

The Jacobian of~\eqref{eq:odes} is given by
\[
J=\begin{pmatrix}
	2x & 2\alpha \\ 2\beta & 2y
\end{pmatrix}.
\]
Hopf bifurcations occur when $\Tr (J)=0$, $\det(J)>0$; there is a neutral saddle if $\Tr (J)=0$ and $\det(J)<0$. 
First, we compute the parameter values required to have $\Tr (J)=0$. This requires $x=-y$ at the equilibrium. It is straightforward to show that this occurs when $x=-\gamma/(\alpha+\beta)$, and so
\[
\lambda=\lambda_{{\Tr}0}
\]
where $\lambda_{{\Tr}0}$ is given in (\ref{eq:lambdaTr0}). In the mutualistic case $\alpha\beta>0$, $x=-y$ implies that $\det(J)=-4(x^2+\alpha\beta)<0$. In this case Hopf bifurcations cannot occur. In the mixed case $\alpha\beta<0$, Hopf bifurcations will occur when $x^2>-\alpha\beta$, that is, when $\gamma^2>-\alpha\beta(\alpha+\beta)^2$.

At a Takens--Bogdanov bifurcation, $\Tr (J)=\det(J)=0$. This requires in addition that $xy=\alpha\beta$. Simple algebraic manipulations show that in parameter space this occurs when
\begin{align}
    \label{eq:TB}
\gamma=\gamma_{TB\pm}\equiv \pm\sqrt{-\alpha\beta}(\alpha+\beta),
\end{align}
which is, as expected, at the boundary of the existence conditions for the Hopf bifurcations, and requires $\alpha\beta<0$.
Combining this information shows that there is a Hopf bifurcation along the curve $\lambda=\lambda_{{\Tr}0}$ when $\gamma_{TB-}<\gamma<\gamma_{TB+}$. If $\gamma$ is outside this range then the curve $\lambda=\lambda_{{\Tr}0}$ corresponds to a neutral saddle.

The criticality of the Hopf bifurcation can be computed using the first Lyapunov coefficient $\ell_1$ \cite{kuznetsov2004elements}, which equals
\begin{align}\label{eq:hopf:l1}
    \ell_1&=\frac{4\beta\gamma}{\omega^3}\mbox{,\quad where\quad}\omega=\sqrt{\det(J)}=2\sqrt{-\gamma^2-\alpha\beta(\alpha+\beta)^2},
\end{align} 
for the parameter range $\alpha>|\beta|>0>\beta$  we consider (see \cite{SNinteraction-2025-pages} for computer algebra script). The expression for $\ell_1$ depends on the scaling of the eigenvectors of the linearization, but its sign is independent of this scaling. We used the eigenvector basis $((\omega,0)^\tran,J(\omega,0)^\tran)$ for \eqref{eq:hopf:l1}, where $\omega$ is the imaginary part of the critical eigenvalues at the Hopf bifurcation. 

Thus, the Hopf bifurcation is supercritical for $\gamma>0$ and subcritical for $\gamma<0$. 
At $\gamma=0$ there is a Bautin point (the first Lyapunov coefficient $\ell_1$ is equal to zero). Here, the second Lyapunov coefficient can be computed, and this will be non-zero unless $\beta=\pm\alpha$.

For our bifurcation diagrams, we will use  a new parameter $\mu$, which measures the distance of $\lambda$ from its value when the  trace equals zero: 
\begin{equation}
    \mu=\lambda-\lambda_{{\Tr}0}=\lambda-\frac{\gamma(\gamma+\alpha^2-\beta^2)}{(\alpha+\beta)^2}
    \label{eq:mu_def}
\end{equation} 
where $\lambda_{{\Tr}0}$ is defined in (\ref{eq:lambdaTr0}).
In these new parameters $\mu,\gamma,\alpha,\beta$ there is a Hopf bifurcation at $\mu=0$ for any $\alpha\beta<0$ and any $\gamma_{TB-}<\gamma<\gamma_{TB+}$.

\subsection{Saddle-node and cusp bifurcations}

The locus of saddle-node bifurcations in the parameter space system~\eqref{eq:odes} can be parametrised by $x$ as
\begin{align*}
    \lambda&=\frac{x^2}{2}+\beta x+\frac{\alpha^2\beta}{x}+\frac{\alpha^2\beta^2}{2x^2},&
    \gamma&=\frac{x^2}{2}-\beta x+\frac{\alpha^2\beta}{x}-\frac{\alpha^2\beta^2}{2x^2},&y&=\frac{\alpha\beta}{x},
\end{align*}
resulting in two disjoint curves, one for $x>0$ and one for $x<0$ (see Figure~\ref{fig:bifset-abp}).

It is also possible to compute the locus in $(\alpha,\beta,\gamma)$ such that there is a cusp point for some value of the parameter $\mu$, or equivalently $\lambda$ (the cusp is called a hysteresis point in \cite{golubitsky1985}). We refer to \cite{golubitsky1985} for details (although there is a sign error in their expression that we correct here). In summary, we find there is a cusp point precisely when the parameters are on the following locus:
\begin{align}\label{eq:cusp}\lambda&=\lambda_{C}\equiv\frac{3}{2}\alpha^{\frac{2}{3}}\beta^{\frac{2}{3}}(\alpha^{\frac{2}{3}}+\beta^{\frac{2}{3}}),&
\gamma&=\gamma_{C}\equiv\frac{3}{2}\alpha^{\frac{2}{3}}\beta^{\frac{2}{3}}(\alpha^{\frac{2}{3}}-\beta^{\frac{2}{3}}).
\end{align}
\citet{golubitsky1985} also give formulae for double limit points and for bifurcations of transcritical type; however, these are dependent on the choice of distinguished parameter.

\section{Dynamics of the mutualistic case}
\label{sec:mutual}

\begin{figure}
	\setlength{\unitlength}{1mm}
    \begin{center}
    \begin{picture}(160,75)
 \put(0,0){\includegraphics[height=7.7cm]{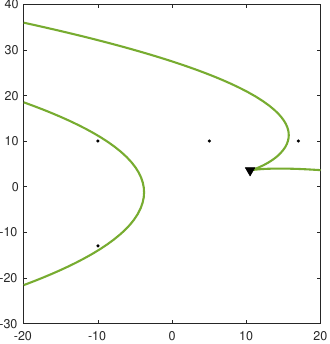}}
 \put(84,3.5){\includegraphics[height=7.3cm]{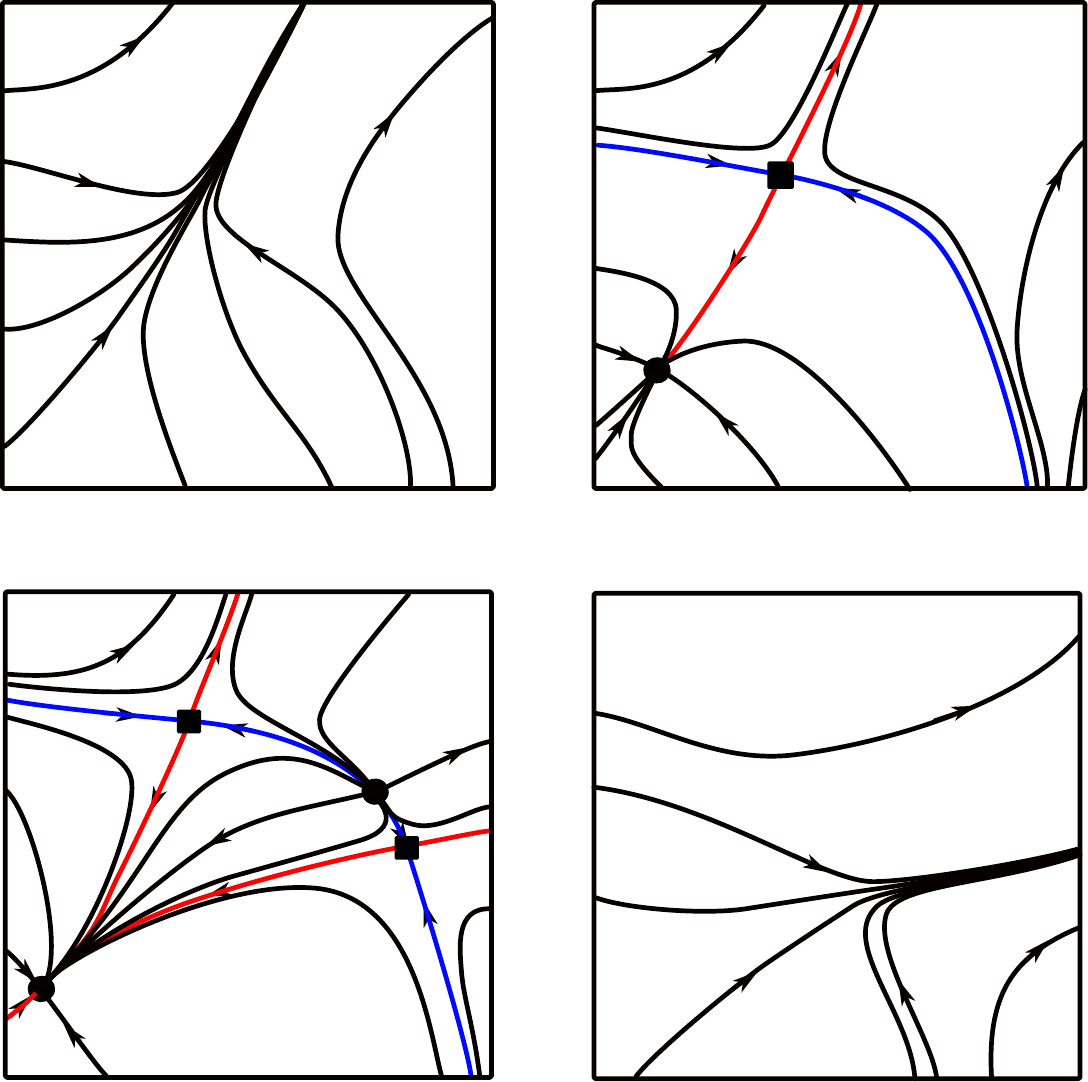}}
    \put(67,-1){$\mu$}
    \put(-1,70){$\gamma$}
    \put(16,43){(a)}
    \put(45,41){(b)}
    \put(64,41){(c)}
    \put(16,23){(d)}

    \put(78,72){(a)}
    \put(118,72){(b)}
    \put(78,33){(c)}
    \put(118,33){(d)}
    
   \end{picture}
    \end{center}
    \caption{In the left panel, the green curves show locations of saddle-node bifurcations in $\mu$-$\gamma$ space for system (\ref{eq:odes}) in the mutualistic case ($\alpha=3.1$, $\beta=1.3$). The triangle indicates a cusp point. Panels (a)-(d) show phase portrait sketches for  $(x,y)\in[-10,10]^2$  at the labelled parameter values in the left panel. There are no equilibria in (a,d), two equilibria in (b) and four equilibria in (c).  Black lines indicate trajectories, squares indicate saddles (where the stable manifold is blue, the unstable manifold is red) and disks indicate sources/sinks. In (a,d), a ghost attractor exists that passes through the region close to a saddle-node bifurcation. In (b,c), there are trapping regions consisting of the local basin of attraction of the stable node in the lower left.}
    \label{fig:bifset-abp}
\end{figure}
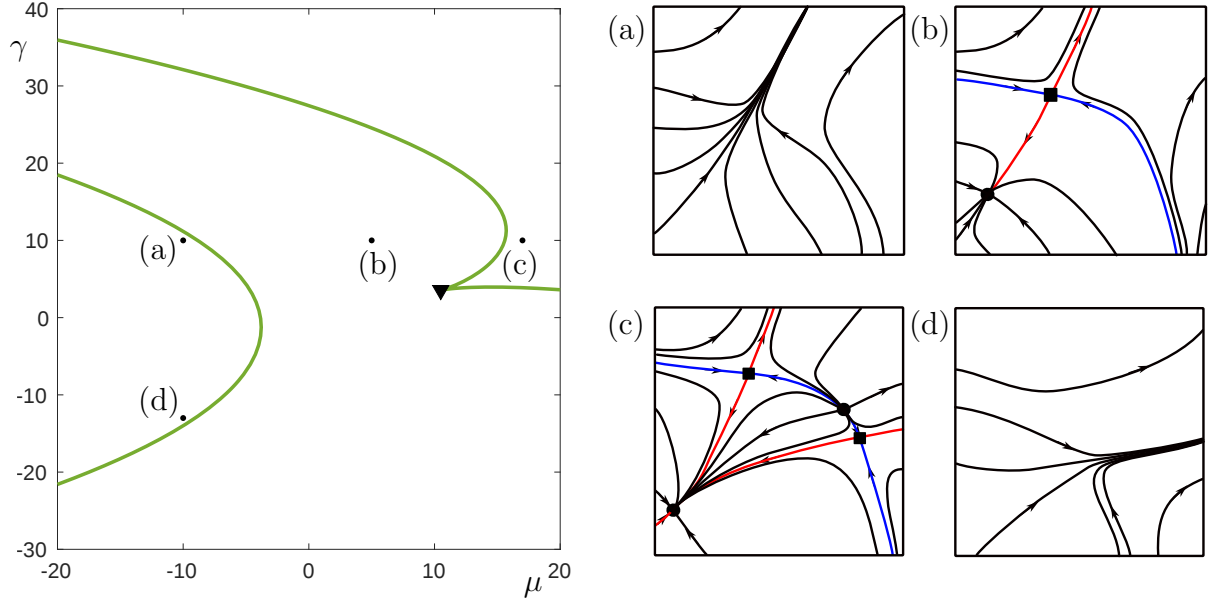

The mutualistic case has $\alpha\beta>0$. 
We show in Figure~\ref{fig:bifset-abp}  a typical arrangement of the curves of saddle-node bifurcations when $\alpha\beta>0$. These curves divide parameter space into three regions. Representative phase portraits for the indicated points are also shown.

For the region with points (a,d) there are no equilibria. However, since the points (a,d) are very close to a saddle-node bifurcation in the parameter plane,  we can see that many trajectories with initial conditions in the lower left exit the region of parameter space shown at a very similar position; this is because all of these trajectories are funnelled through the same region of phase space. At points (b) and (c) (and throughout the regions containing these points) there is a stable equilibrium such that some trajectories are trapped in its basin of attraction.

\section{Dynamics of the mixed case}
\label{sec:mixed}

The mixed case has $\alpha\beta<0$. We first give a description of the bifurcation curves in the $\mu$-$\gamma$ plane for fixed $\alpha,\beta$, in Section~\ref{sec:mixed_bif}. Then in Section~\ref{sec:mixedpp} we describe the dynamics that occur in each region of parameter space.

\subsection{Bifurcation curves}
\label{sec:mixed_bif}

\begin{figure}
	\setlength{\unitlength}{1mm}
    \begin{center}
    \begin{picture}(146,70)
     \put(84,1){\includegraphics[height=3.8cm]{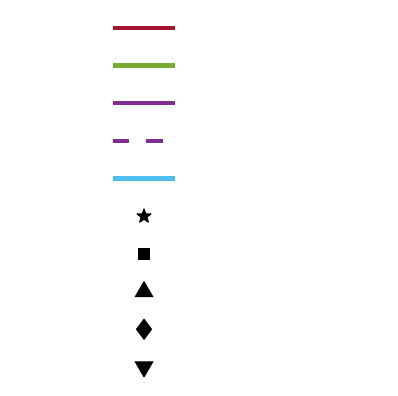}}
        \put(0,0){    \includegraphics[height=7cm]{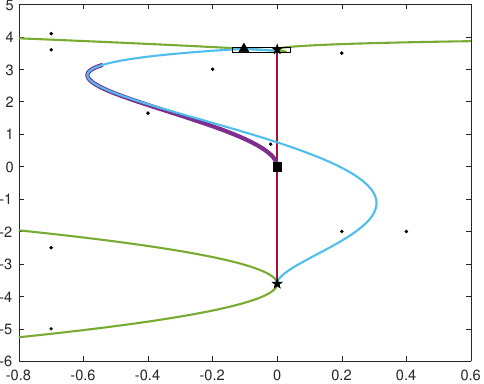}}
      \put(91,39){    \includegraphics[height=3cm]{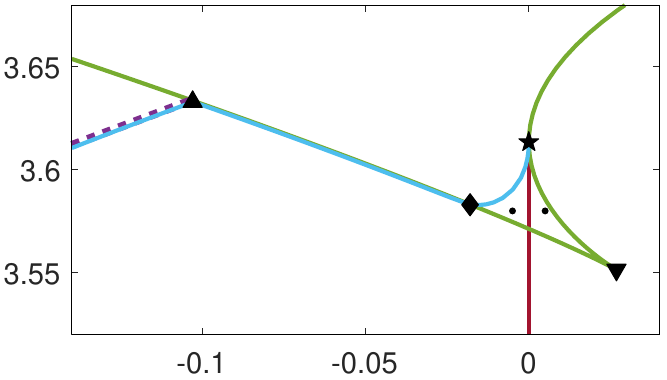}}
    \put(82,-1){$\mu$}
    \put(-1,66){$\gamma$}
    \put(143,39){$\mu$}
    \put(94,66.5){$\gamma$}
    \put(11,65){(a)}
    \put(11,58){(b)}
    \put(11,22){(c)}
    \put(11,11){(d)} 
    \put(40,53){(f)}
    \put(44,40){(g)}
    \put(59,56){(h)}
    \put(59,30){(i)}
    \put(75,30){(j)}
    \put(22,45){(e)}
   \put(128.5,47){(k)}
   \put(139,51){(l)}
   
   \put(100,55.5){\scriptsize (1)}
   \put(106,58){\scriptsize (2)}
   \put(118,58){\scriptsize (3)}
   \put(125,50.5){\scriptsize (4)}
   \put(129.5,55){\scriptsize (5)}

    \put(104,35.7){\scriptsize Hopf}
    \put(104,32){\scriptsize Saddle-node of equilibria}
   \put(104,28.2){\scriptsize Saddle-node of periodic orbits (SNPO)}
   \put(104,24.7){\scriptsize SNPO (conjectured)}
   \put(104,21){\scriptsize Long-period periodic orbits}
   \put(104,17.3){\scriptsize Takens--Bogdanov point}
   \put(104,13.8){\scriptsize Bautin point}
   \put(104,10.4){\scriptsize SNICeroclinic}
   \put(104,6.7){\scriptsize Non-central SNIC}
   \put(104,2.9){\scriptsize Cusp point}
    \end{picture}
    \end{center}
    \caption{Bifurcation curves in $\mu$-$\gamma$ space when $\alpha=3.1$, $\beta=-1.3$. The upper right-hand panel shows a zoom of the enclosed area in the left-hand panel. Curves of codimension one bifurcations are coloured as in the legend.
    The curve of long-period periodic orbits approximates a homoclinic orbit everywhere except where it coincides with the curve of saddle-node bifurcations, where there is a SNIC bifurcation. 
    Symbols indicate codimension-2 points. The curve of SNPO is conjectured to terminate in the SNICeroclinic, the conjectured curve is shown as a dashed purple line in the right-hand panel (details in the text). The labels with numbers in the zoom indicate the schematic phase portraits in Figure~\ref{fig:SNICs}. The dots labelled with letters correspond to parameter values for phase portraits in Figure~\ref{fig:pp-abn}. Reproducing code and interactive version of figure are available at \cite{SNinteraction-2025-pages}.}
    \label{fig:bifset-abn}
\end{figure}

When $\alpha\beta<0$ then the set of bifurcations is considerably more complicated than in the mutualistic case. 
In addition to the curves of saddle-node bifurcation, recall from Section~\ref{sec:hopf} that there is a curve of Hopf bifurcations for $\mu=0$, $\gamma_{TB-}<\gamma<\gamma_{TB+}$. At each end of this curve, there is a Takens--Bogdanov (TB) point. At $\gamma=0$, the Hopf bifurcation has zero first Lyapunov coefficient; this is a Bautin point. The existence of these codimension-two TB and Bautin points means there must exist curves of homoclinic bifurcations and saddle-nodes of periodic orbits.
 These curves cannot be computed analytically, so we use the continuation software COCO~\cite{dankowicz2013recipes} to follow these branches numerically. The curve of homoclinic orbits is approximated by following a branch of periodic orbits with fixed large period ($T\approx170$). A typical set of bifurcation curves is shown in Figure~\ref{fig:bifset-abn} and corresponding phase portraits are shown in Figure~\ref{fig:pp-abn}.

At each TB point, a curve of saddle-node bifurcations meets a curve of Hopf bifurcations, and a curve of homoclinic bifurcations emanates from each of the two TB points.       
These two curves of homoclinic bifurcations form a single continuous curve
between the two TB points such that the long-period approximation is a single family or periodic orbits. As shown in the zoomed-in panel of Figure~\ref{fig:bifset-abn}, the curve of homoclinics collides with the curve of saddle-node bifurcations and turns into a short branch of SNIC bifurcations. The branch of SNIC bifurcations is bounded on the right in a non-central SNIC, and on the left in a SNICeroclinic bifurcation~\cite{nechyporenko2024}.

\begin{figure}
	\setlength{\unitlength}{1mm}
    \begin{center}
    \begin{picture}(160,40)
     \put(0,0){\includegraphics[width=\textwidth]{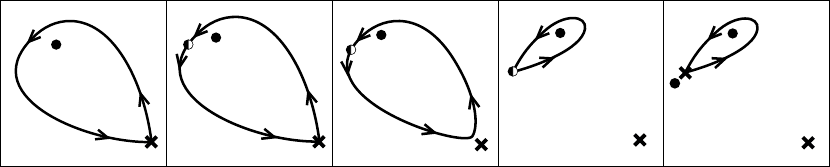}}
    \put(1,29){(1)}
    \put(34,29){(2)}
    \put(67,29){(3)}
    \put(100,29){(4)}
    \put(133,29){(5)}
    \put(2,2){\vector(1,0){5}}
    \put(2,2){\vector(0,1){5}}
    \put(1.6,8.2){$y$}
    \put(7.5,1.6){$x$}
    \end{picture}
    \end{center}
    \caption{Each panel shows a schematic of the phase portrait along the homoclinic and SNIC curves (light blue) shown in the right-hand panel in Figure~\ref{fig:bifset-abn}. Equilibria are shown by various symbols: crosses are saddles, black circles are nodes, half-black circles are saddle-nodes. (1) on the homoclinic bifurcation curve to the left of the SNICeroclinic bifurcation; (2) at the SNICeroclinic bifurcation; (3) along the SNIC bifurcation; (4) at the non-central SNIC bifurcation; (5) on the homoclinic bifurcation curve to the right of the non-central SNIC.}
    \label{fig:SNICs}
\end{figure}
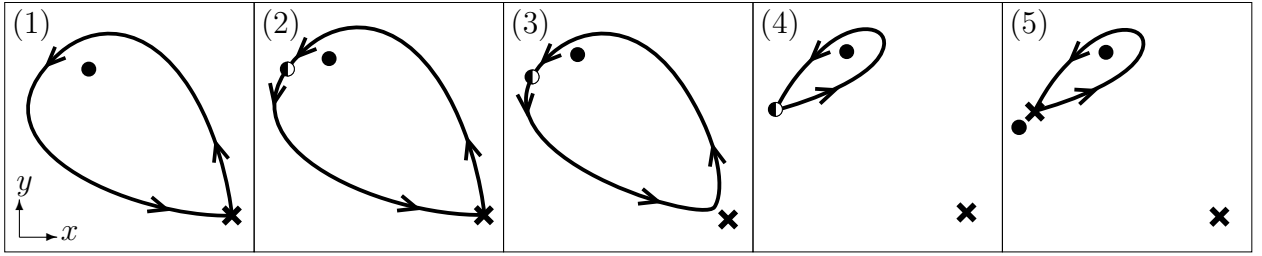

In Figure~\ref{fig:SNICs} we show schematic phase portraits of the homoclinic and SNIC orbits as the curve of homoclinic bifurcations passes through the SNICeroclinic and non-central SNIC transition points in the parameter plane. The five points are indicated by (1)-(5) along the light-blue curve in Figure~\ref{fig:bifset-abn}. We note that the homoclinic orbits to the left and right of the SNIC curve are homoclinic to different saddle points; the nearby periodic orbits also have different stabilities. 
The homoclinic curve to the left of the SNICeroclinic is shown in panel (1); it is a large homoclinic orbit to a saddle in the lower right of the phase portrait. At the SNICeroclinic bifurcation, this homoclinic orbit collides with a saddle-node equilibrium (panel (2)). 
On the SNIC curve, the orbit is detached from the saddle in the lower right (panel (3)), and as we move to the right along this curve, the orbit changes shape and becomes smaller. 
At the non-central SNIC, the orbit again looks like a homoclinic orbit (panel (4)) to a saddle-point, although here the equilibrium is a non-hyperbolic saddle-node. The homoclinic orbits on the right of the SNIC curve (panel (5)) are homoclinic to a saddle on the left of the phase portrait. These homoclinic bifurcations disappear in the nearby Takens--Bogdanov bifurcation.

SNICeroclinic bifurcations have recently been studied in~\cite{nechyporenko2024}. However, in that work, they only considered the case when both the SNIC and homoclinic bifurcations result in periodic orbits of the same stability. In our case, the homoclinic bifurcation to the left of the SNIC is subcritical, i.e. it 
produces an unstable periodic orbit.
The SNIC bifurcation, however, produces a stable periodic orbit. Thus, there must also be a curve of SNPOs emanating from the SNICeroclinic point. We conjecture that this is the same curve of SNPO as that emanating from the Bautin point. Numerical continuation of SNPOs close to homoclinic bifurcations is challenging due to the difficulties in controlling Floquet multipliers of orbits which spend a long time near equilibrium solutions. Due to this, we were unable to continue the curve of SNPOs to the SNICeroclinic point, but we conjecture that this is where it must terminate (dashed purple curve in inset to Figure~\ref{fig:bifset-abn}.)

\subsection{Phase plane dynamics}
\label{sec:mixedpp}

\begin{figure}
\setlength{\unitlength}{1mm}
    \centering
   \begin{picture}(170,110)
    \put(6,0){\includegraphics[width=158mm]{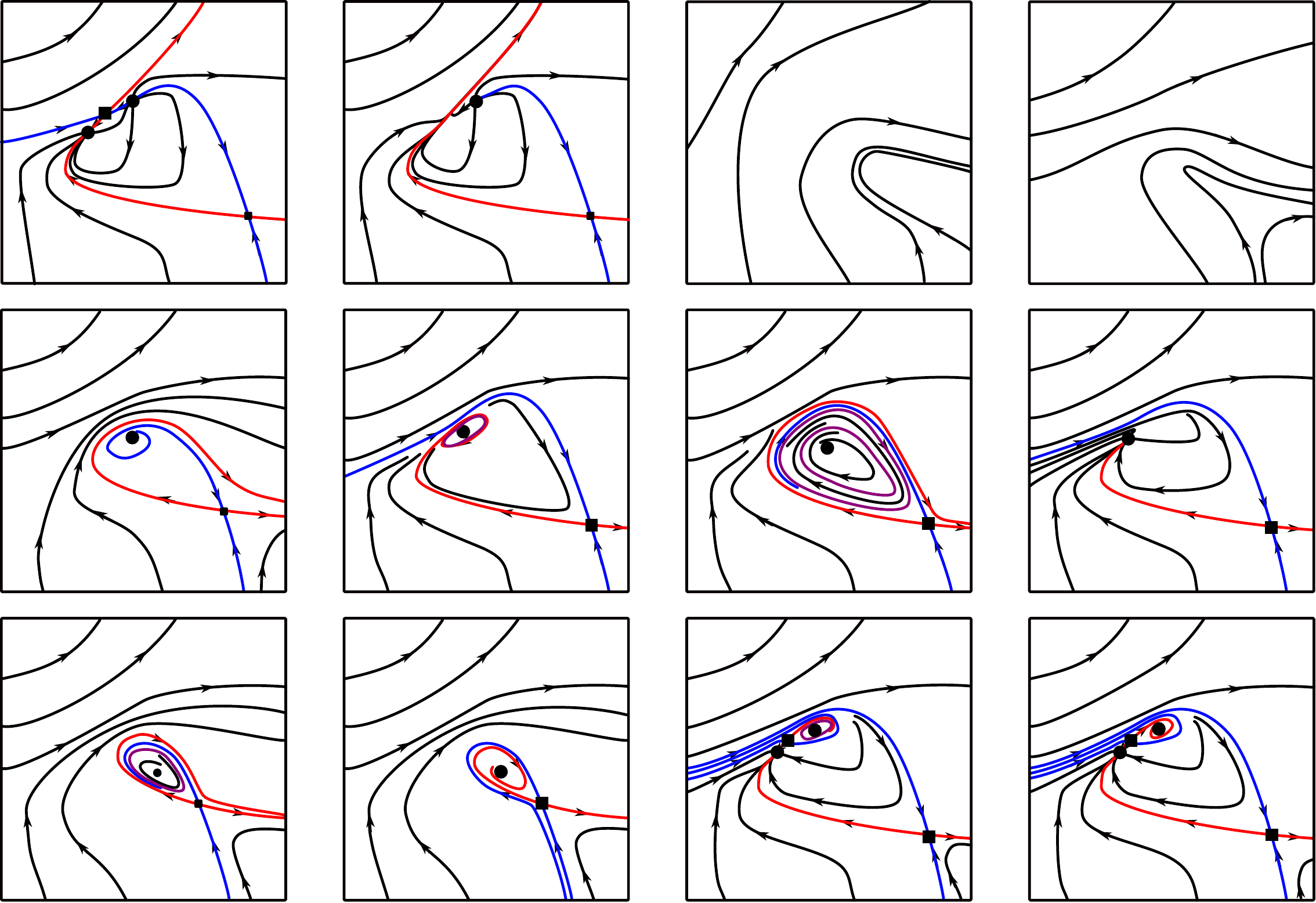}}
    \put(0.5,30){(i)}
    \put(42,30){(j)}
    \put(83,30){(k)}
    \put(124.5,30){(l)}
      \put(0,67){(e)}
    \put(41.5,67){(f)}
    \put(83,67){(g)}
    \put(124,67){(h)}
      \put(0,105){(a)}
    \put(41.5,105){(b)}
    \put(83,105){(c)}
    \put(124,105){(d)}
    \end{picture}  
    \caption{Phase portrait sketches of (\ref{eq:odes}) for the parameter values corresponding to the points shown in Figure~\ref{fig:bifset-abn}(a)-(l). The region $(x,y)\in[-10,10]^2$ is shown.  The lines are as in Figure~\ref{fig:bifset-abp}; in addition periodic orbits in (f), (g), (i) and (k) are shown in purple. Note that (e) has an unstable periodic orbit that encloses a stable periodic orbit. These merge at a saddle-node bifurcation of period orbits on passing from (g) to (e). }
    \label{fig:pp-abn}
\end{figure}
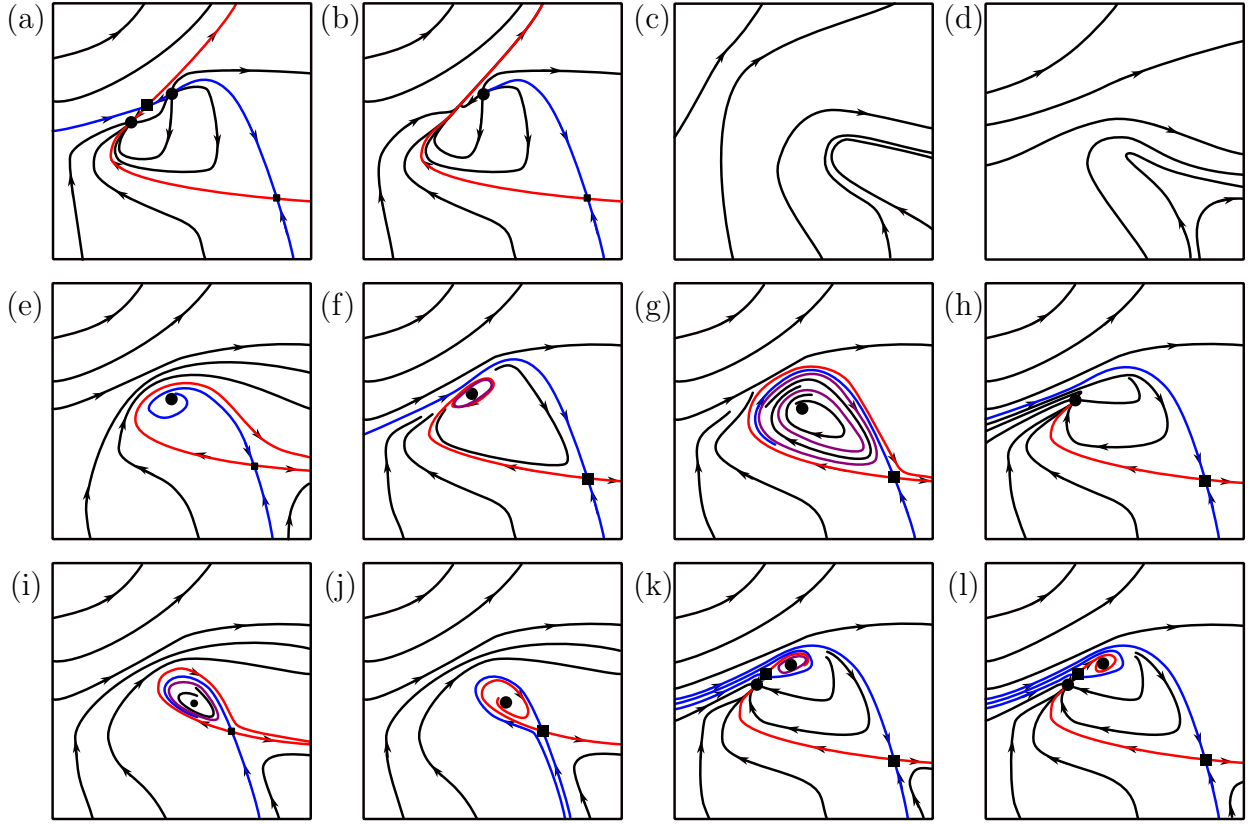

In Figure~\ref{fig:pp-abn}(a)-(l) we show the phase plane dynamics at the points indicated in Figure~\ref{fig:bifset-abn}.

At point (c), there are no equilibrium points. If $\gamma$ is decreased to point (d) close to the saddle-node bifurcation, there is a large set of trajectories that exit the box near the lower right. This is indicative of a ghost attractor~\cite{koch2024ghost} (a ``nearly invariant'' region of the phase space, where trajectories spend a long time before leaving). As $\gamma$ is increased from point (c) and the saddle-node bifurcation is crossed, two equilibria are created, but both are unstable. At point (b), we are now close to the second saddle-node bifurcation. This creates a ghost attractor as in (d). 
As with Figure~\ref{fig:bifset-abp}(a), the ghost attractors in (d) and (b) are identified by the large set of trajectories which have a very similar exit point from the region of phase space shown in the figure. This is discussed more in Section~\ref{sec:transient}.

From the region containing point (b), as $\mu$ is increased, periodic orbits are created and destroyed in multiple different bifurcation types. We begin with those that occur for $\gamma<0$. 
As $\mu$ is increased past $0$ to point (i), one of the equilibria gains stability and an unstable periodic orbit is created in a subcritical Hopf bifurcation. Trajectories that start inside the periodic orbit will be trapped. As $\mu$ is increased further to point (j), the periodic orbit disappears in a subcritical homoclinic bifurcation; the stable equilibrium now has a larger basin of attraction. Point (h) is in the same region, but the basin of attraction of the stable equilibrium is even larger.

For $\mu>0$, the Hopf bifurcation is supercritical. Thus at point (g) there are two periodic orbits, a stable one inside an unstable one. The unstable periodic orbit again forms the boundary of an inaccessible region of phase space. These periodic orbits collide and disappear as the saddle-node of periodic orbits (SNPO) bifurcation is crossed to point (e). At point (e) there is a ghost attractor caused by the nearby SNPO; trajectories inside where the periodic orbits used to be exit the region of phase space at a similar point. If instead we increase $\gamma$ from point (g) to point (f), we cross the subcritical homoclinic bifurcation, and the unstable periodic orbit disappears while the stable periodic orbit remains. The boundary of the basin of attraction of this periodic orbit comprises the stable manifold of the saddle equilibrium. Many trajectories are now unable to escape, as they are in the basin of attraction of the stable periodic orbit.

We can move from point (f) to point (a) by increasing $\gamma$ across the SNIC bifurcation; the stable periodic orbit disappears, but the overall escape dynamics are largely unchanged. Points (k) and (l) similarly have very similar escape dynamics. In (k), the stable periodic orbit remains; a stable equilibrium has been created in the saddle-node bifurcation.  As we cross the supercritical Hopf bifurcation to point (l), the stable periodic orbit disappears; now there are two stable equilibria. In both (k) and (l),  the trapping region consists of the union of the basins of attraction of these stable objects.

\subsection{Codimension-two bifurcations}

The final part of our analysis of the mixed case is to identify the location of the codimension-two points within the full parameter space. As described in Section~\ref{sec:scaling}, we can fix one parameter by an appropriate scaling. So we fix $\alpha$ and vary $\beta$, along with $\mu$ and $\gamma$ such that codimension-two bifurcations exist along curves in the three-parameter $(\mu,\beta,\gamma)$-space . 

We already have analytical expressions for the cusp bifurcations \eqref{eq:cusp}, Takens--Bogdanov bifurcations \eqref{eq:TB}, and Bautin points (see Section~\ref{sec:hopf}). We find the curves of SNICeroclinic and non-central SNIC points using continuation in COCO.

\paragraph{SNICeroclinic and non-central SNIC}

The SNICeroclinic and non-central SNIC points are found and tracked numerically. The SNICeroclinic is approximated by an orbit segment of large fixed integration period that starts at $u_\mathrm{sn}+s_0v_\mathrm{c}$, and ends at $u_\mathrm{sa}+s_1v_\mathrm{s}$. The point $u_\mathrm{sn}=(x_\mathrm{sn},y_\mathrm{sn})$ is the saddle-node equilibrium, $v_\mathrm{c}$ is the eigenvector for eigenvalue $0$ of the linearization in $u_\mathrm{sn}$. The point $u_\mathrm{sa}=(x_\mathrm{sa},y_\mathrm{sa})$ is the saddle equilibrium, $v_\mathrm{s}$ is the eigenvector for the negative eigenvalue of the linearization in $u_\mathrm{sa}$. The numbers $s_0$ and $s_1$ are small but variable along the curve in $(\mu,\beta,\gamma)$-space.

The non-central SNIC is approximated by an orbit segment $u([0,T])=(x([0,T]),y([0,T]))$
 of large fixed integration period $T$ that starts at $u_\mathrm{sn}+s_0v_\mathrm{c}$ (so, $0=u_\mathrm{sn}+s_0v_\mathrm{c}-u(0)$), and ends near $u_\mathrm{sn}$ satisfying $w_\mathrm{c}^\tran[u_\mathrm{sn}-u(T)]=0$. The point $u_\mathrm{sn}=(x_\mathrm{sn},y_\mathrm{sn})$ is the saddle-node equilibrium, $v_\mathrm{c}$ is the (right) eigenvector for eigenvalue $0$ of the linearization in $u_\mathrm{sn}$, and $w_\mathrm{c}$ is the corresponding left eigenvector.
 The details of the defining systems for the numerical approximation of non-central SNIC and SNICeroclinic are given in Appendix~\ref{sec:defining:systems}. 
 
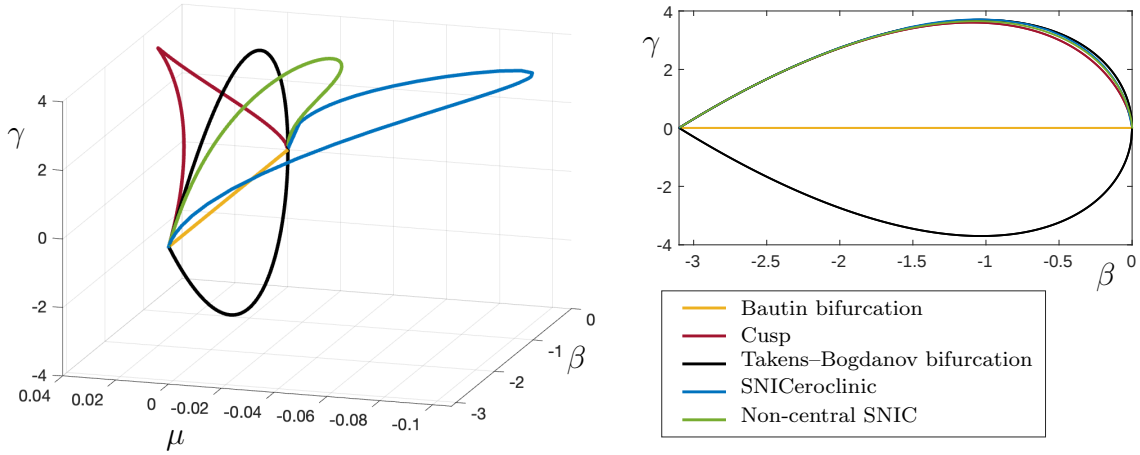
\begin{figure}
	\setlength{\unitlength}{1mm}
    \begin{center}
    \begin{picture}(160,52)

\put(0,0){\includegraphics[height=6.5cm]{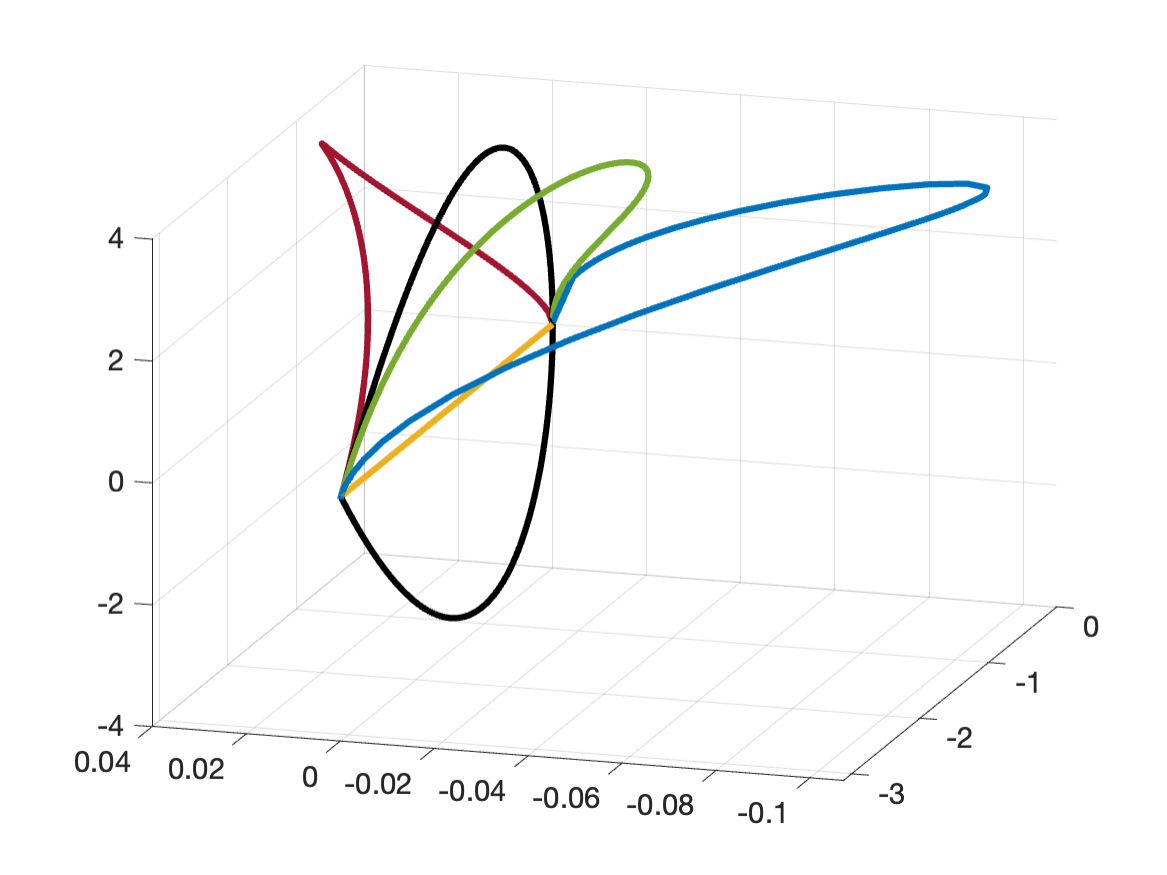}}
\put(78,12){$\beta$}
\put(4,42){$\gamma$}
 \put(25,2){$\mu$}

 \put(90,25){\includegraphics[height=3.5cm]{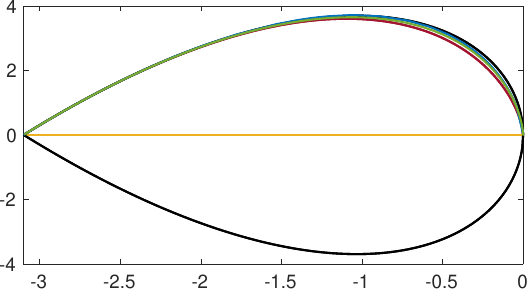}}
\put(148,23){$\beta$}
\put(88,54){$\gamma$}

 \put(90.5,2.3){\includegraphics[height=2cm]{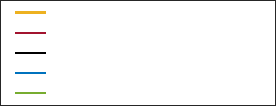}}
 \put(101,19.0){\scriptsize Bautin bifurcation}
   \put(101,15.5){\scriptsize Cusp}
   \put(101,12.3){\scriptsize Takens--Bogdanov bifurcation}
   \put(101,8.7){\scriptsize SNICeroclinic}
\put(101,4.9){\scriptsize Non-central SNIC}
    
   \end{picture}
    \end{center}
    \caption{Curves of codimension-2 bifurcations. The left panel shows the curves in $\mu$-$\beta$-$\gamma$, and the right panel shows the projections onto the $\beta$-$\gamma$ plane. Reproducing code and interactive version of figure are available at \cite{SNinteraction-2025-pages}.
    \label{fig:cod2-bifs}}
\end{figure}

The resulting approximations for SNICeroclinic (in blue) and non-central SNIC (in green) are shown in Figure~\ref{fig:cod2-bifs} as curves in $(\mu,\beta,\gamma)$-space, alongside the curves for the cusp bifurcations, Takens--Bogdanov bifurcations and Bautin points.

\paragraph{Implications for codimension-one bifurcations} Figure~\ref{fig:cod2-bifs} does not show the codimension-one bifurcations: these form surfaces described below that join the curves of codimension-two points. The surface of Hopf bifurcations is the part of the plane $\mu=0$ enclosed by the closed curve of Takens-Bogdanov bifurcations. The line of Bautin bifurcations at $\mu=\gamma=0$ is inside this Hopf bifurcation surface. The surface of saddle-node bifurcations of equilibria passes through the SNICeroclinic  (blue curve), the Non-central SNIC (green curve), the Takens-Bogdanov bifurcation (closed black curve) and the cusp (dark red curve). Its surface tangent is parallel to the surface of Hopf bifurcations at the plane $\mu=0$ in the Taken-Bogdanov bifurcations. The surface of homoclinic bifurcations (homoclinic to a saddle equilibrium) has boundaries on the Takens-Bogdanov curve, the Non-central SNIC curve and the SNICeroclinic curve. The surface of SNIC bifurcations is part of the surface of saddle-node bifurcations of equilibria, bounded by the Non-central SNIC and SNICeroclinic bifurcation.

\section{Implications for transient behaviour}
\label{sec:transient}

In the context of coupled saddle-nodes with a global reinjection mechanism (as considered by Baesens and MacKay~\cite{baesens2013interaction}), we are concerned as to whether all trajectories get through the region of phase space where the saddle-node bifurcations will take place. As can be seen in the phase space diagrams, and we will detail further in this section, there are some parameter regions where all trajectories get through eventually, others where all trajectories are trapped, and a third class where a subset of trajectories are trapped, and others can pass through. 

Near the boundaries between parameter sets where trajectories are trapped, and those where trajectories are able to pass through, we see ghost attractors~\cite{deco2012ongoing,koch2024ghost,morozov2024long}, also known as bottlenecks~\cite{strogatz2024nonlinear}. 
In this section, we give some numerical results and interpret these in the context of ghost attractors. We note that, particularly in the mixed case, being able to correctly identify the boundaries of the regions of parameter space where all trajectories `get through' requires knowing about the global bifurcations investigated in Section~\ref{sec:mixed} above; in particular, knowing only the existence criteria of the equilibria (as in~\cite{golubitsky1985}) is not sufficient in the mixed case as a variety of bifurcations involving periodic orbits are possible.

In order to visualise the presence or absence of ghost attractors, for each parameter set, we compute trajectories for $150^2$ evenly spaced initial conditions in the box $(x,y)\in [-10,10]^2$. For each trajectory, we calculate the time taken to leave this box, and additionally, the point on the boundary through which the trajectory leaves the box. In all cases, trajectories leave either through the top or the right-hand side of the box. We use this exit position to define an `departure angle' for each trajectory, which varies linearly from $-1$ to $0$ from left to right along the top of the box, and from $0$ to $+1$ from top to bottom along the right-hand side of the box. More precisely, we measure the angle $\psi$ from the bottom left corner to the exit point of each trajectory on the right or top side of the box and characterise each trajectory that exits these sides with the scalar quantity $4\psi/\pi-1$ which is in the interval $[-1,1]$.

\subsection{Transients in the mutualistic case}

\begin{figure}[ht!]
\setlength{\unitlength}{1mm}
    \begin{center}
    \begin{picture}(165,102)(0,0)

\put(4,48){\includegraphics[width=0.47\linewidth]{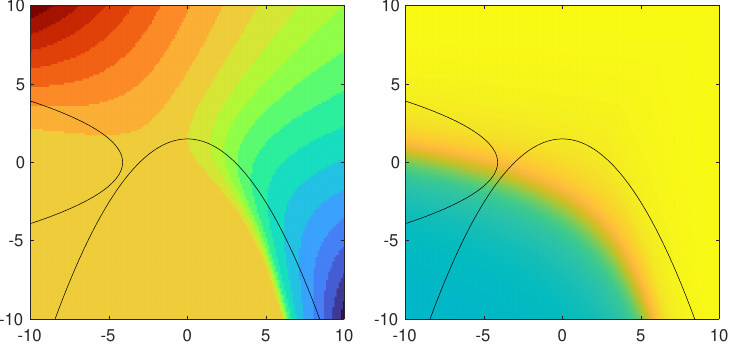}}
\put(86,48){\includegraphics[width=0.47\linewidth]{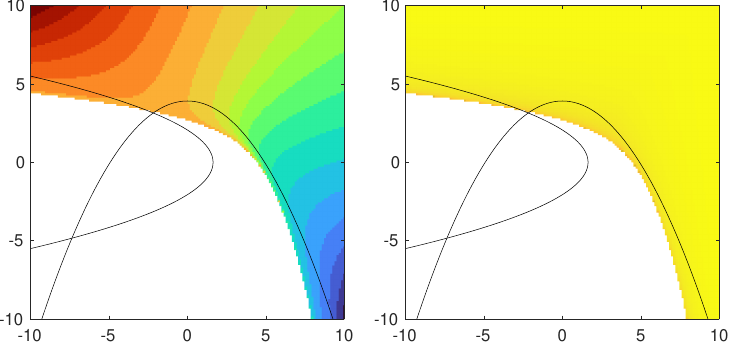}}
\put(4,2){\includegraphics[width=0.47\linewidth]{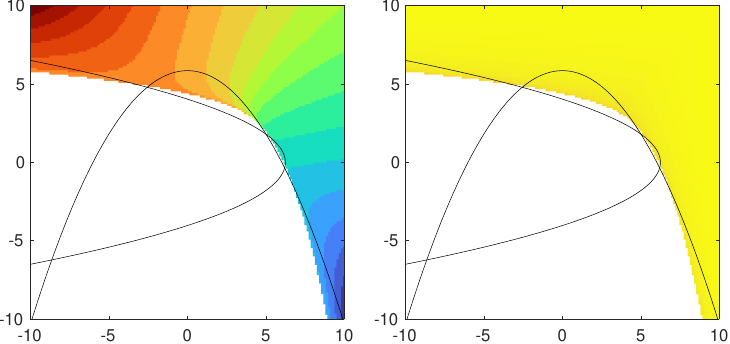}}
\put(86,2){\includegraphics[width=0.47\linewidth]{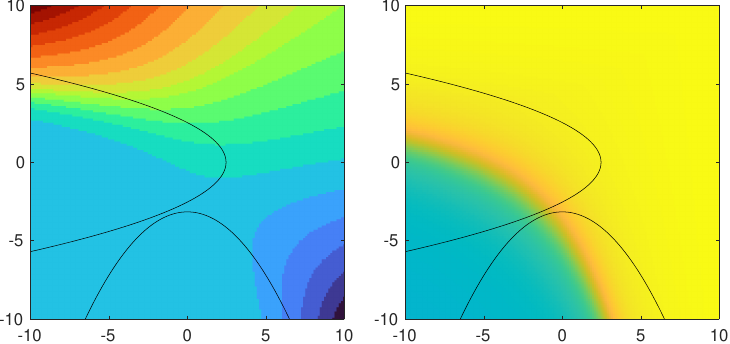}}

 \put(0,88){(a)}
 \put(37,46){$x$} 
\put(2,83){$y$}
\put(77,46){$x$} 
 
 \put(82,88){(b)}
\put(119,46){$x$} 
\put(84,83){$y$}
\put(159,46){$x$} 
 
 \put(0,42){(c)}
\put(37,0){$x$} 
\put(2,35){$y$}
\put(77,0){$x$}
 
 \put(82,42){(d)}
\put(119,0){$x$} 
\put(84,35){$y$}
\put(159,0){$x$} 

\put(37,0){$x$} 
\put(2,35){$y$}
\put(77,0){$x$} 

\put(7,91){\includegraphics[width=0.203\linewidth]{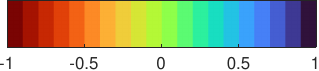}}
\put(10,101){Departure angle}

\put(89,91){\includegraphics[width=0.203\linewidth]{fig-colorbar2.pdf}}
\put(92,101){Departure angle}

\put(47,91){\includegraphics[width=0.203\linewidth]{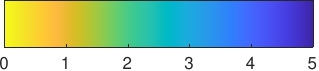}}
\put(53,101){Escape time}

\put(129,91){\includegraphics[width=0.203\linewidth]{fig-colorbar3.pdf}}
\put(135,101){Escape time}
 
 \end{picture}
    \end{center}
    \caption{Each pair of panels shows the phase plane for the mutualistic case, for the parameters indicated in Figure~\ref{fig:bifset-abp}. The black lines in each panel are the nullclines. 
    We start a trajectory at each of a fine grid of points in the box.
    The {\it Departure angle} column indicates the position the trajectory exits the box (as described further in the text). The {\it Escape time} column indicates the time taken for the trajectory to leave (in arbitrary units). Points coloured white correspond to trajectories that never leave the box. In panels (a) and (d), the nullclines do not intersect, and there is a large region of initial conditions which both take a long time to leave the box and have a similar departure angle. In panels (b) and (c), the nullclines intersect, creating a stable equilibrium. The white region corresponds to the basin of attraction of this equilibrium.
    }
    \label{fig:bnfigs-abp}
\end{figure}

In Figure~\ref{fig:bnfigs-abp}(a)-(d), we show the dynamics in the $x$-$y$ plane at each of the labelled points in Figure~\ref{fig:bifset-abp}.  
In panel (a), parameters are such that we are close to a saddle-node bifurcation: this can be seen by the fact that the nullclines are close but do not intersect. The ghost attractor can be clearly seen in the {\it Departure angle}: we can see that a large region of initial conditions in the lower left (orange region) all exit the box with a very similar departure angle. All these trajectories are funnelled through the region of phase space close to where the saddle-node bifurcation will take place. In the {\it Escape time} panel we can additionally see that these trajectories take significantly longer to escape the box than trajectories starting on the other side of the ghost attractor.

In panels (b) and (c), there is a stable equilibrium, and so there are large regions of phase space for which trajectories do not ever exit the box: they are asymptotic to the stable equilibrium. This appears as a white region in the figures, and is bounded by the stable manifold of other saddle points in the phase space. 

Panel (d) is similar to panel (a), but the direction of escape for the trajectories which pass through the ghost attractor is to the right rather than upwards (blue region). There is in fact a continuous path through parameter space (`inside' the left-hand saddle node bifurcation curve in Figure~\ref{fig:bifset-abp}) that can be traversed while continuously varying this exit angle.

\subsection{Transients in the mixed case}

\begin{figure}[ht!]
\setlength{\unitlength}{1mm}
    \begin{center}
    \begin{picture}(165,102)(0,0)

\put(4,48){\includegraphics[width=0.47\linewidth]{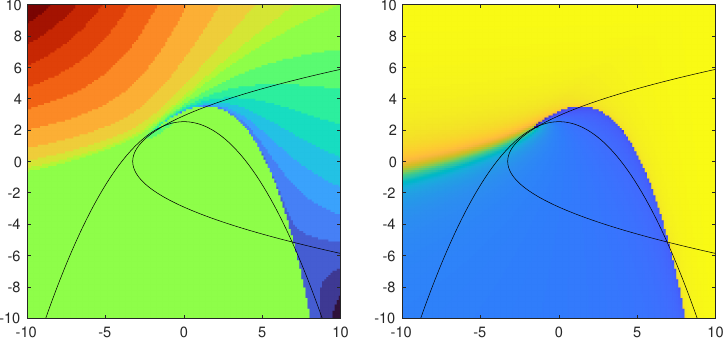}}
\put(86,48){\includegraphics[width=0.47\linewidth]{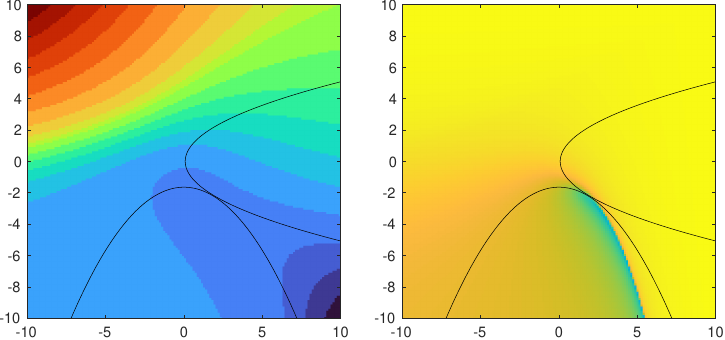}}
\put(4,2){\includegraphics[width=0.47\linewidth]{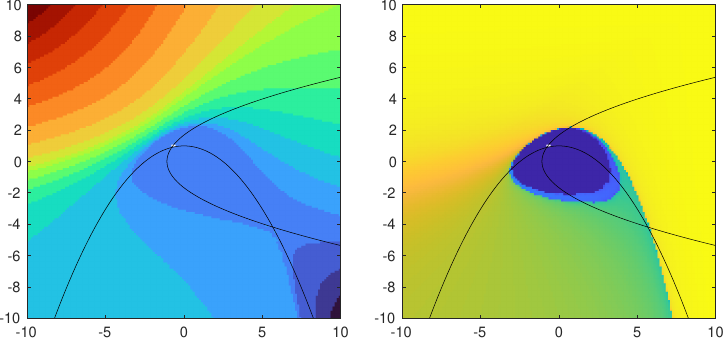}}
\put(86,2){\includegraphics[width=0.47\linewidth]{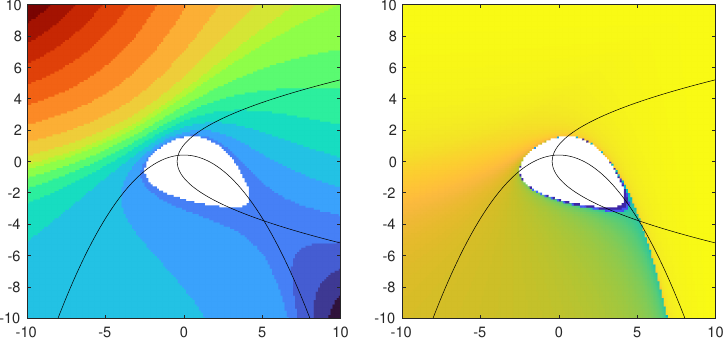}}

 \put(0,88){(b)}
 \put(37,46){$x$} 
\put(2,83){$y$}
\put(77,46){$x$} 
 
 \put(82,88){(d)}
\put(119,46){$x$} 
\put(84,83){$y$}
\put(159,46){$x$} 
 
 \put(0,42){(e)}
\put(37,0){$x$} 
\put(2,35){$y$}
\put(77,0){$x$}
 
 \put(82,42){(g)}
\put(119,0){$x$} 
\put(84,35){$y$}
\put(159,0){$x$}

\put(37,0){$x$} 
\put(2,35){$y$}
\put(77,0){$x$} 

\put(7,91){\includegraphics[width=0.203\linewidth]{fig-colorbar2.pdf}}
\put(10,101){Departure angle}

\put(89,91){\includegraphics[width=0.203\linewidth]{fig-colorbar2.pdf}}
\put(92,101){Departure angle}

\put(47,91){\includegraphics[width=0.203\linewidth]{fig-colorbar3.pdf}}
\put(53,101){Escape time}

\put(129,91){\includegraphics[width=0.203\linewidth]{fig-colorbar3.pdf}}
\put(135,101){Escape time}

 \end{picture}
    \end{center}
    \caption{Each pair of panels shows the phase plane for a subset of the parameters indicated in Figure~\ref{fig:bifset-abn}, namely, the points (b), (d), (e) and (g).  The black lines in each panel are the nullclines. 
    We start a trajectory at each of a fine grid of points in the box.
    The {\it Departure angle} column indicates the position the trajectory exits the box (as described further in the text). The {\it Escape time} column indicates the time taken for the trajectory to leave. Points coloured white correspond to trajectories that never leave the box. (b) and (d) are close to saddle-node bifurcations which results in a larger portion of trajectories being funnelled past a ghost attractor and exiting the box at a similar angle (green and light blue regions respectively). (e) and (g) are either side of a SNPO. 
    In (g), some trajectories are trapped inside a stable periodic orbit. In (e), after the SNPO, all trajectories `inside' the location of the periodic orbits leave the box in a similar direction (medium blue region).
    }
    \label{fig:bnfigs-abn}
\end{figure}

In Figure~\ref{fig:bnfigs-abn}(b), (d), (e) and (g) we show the dynamics in the phase plane for parameter sets corresponding to the matching labelled points in Figure~\ref{fig:bifset-abn}. We choose these points because the phase portraits contain, or are close to containing, ghost attractors. 

In panel (b), we are close to a saddle-node bifurcation. This creates a ghost attractor similar to that seen in Figure~\ref{fig:bnfigs-abp} panels (a) and (d) (green region in the lower left). 
In panel (d), the parameters are also close to a saddle-node bifurcation. The stable equilibrium which is formed in this bifurcation has a fairly small basin of attraction (see panel (j) in Figure~\ref{fig:pp-abn}), and so the region of phase space which passes through the ghost attractor (medium-blue region) is somewhat smaller than in (b). 

In panel (g),  there are two periodic orbits, a stable one inside an unstable one. The unstable periodic orbit again forms the boundary of an inaccessible region of phase space. These periodic orbits collide and disappear as the saddle-node of periodic orbits (SNPO) bifurcation is crossed. Panel (e) has parameters after this SNPO bifurcation. Now all trajectories which start inside the previous location of the periodic orbits take a long time to escape, as they pass through a ghost attractor of the periodic orbits.

\section{Discussion}
\label{sec:discuss}

The richness of bifurcation structure in this simple scenario of two coupled saddle-node normal forms is quite unexpected. We consider this quadratic truncated normal form near the generic situation of an isolated equilibrium with zero eigenvalue and geometric multiplicity two.  The mutualistic case is essentially the local normal form investigated in \cite{baesens2013interaction}. Our analysis of the mixed case is new and contains a number of surprises even at this local level: there are a variety of local bifurcations and bifurcations of connecting orbits that connect up to four equilibria and up to two limit cycles in a number of ways. We have provided analytical formulas for bifurcations of equilibria. Bifurcations of periodic orbits could be determined only using numerical approximation methods; we have developed suitable defining nonlinear boundary-value problems (see Appendix~\ref{sec:defining:systems}) to continue the non-central SNIC and SNICeroclinic bifurcations in normal form coefficients.

The quadratic normal form we consider is sufficient to determine the topological branching in the one-dimensional case. However, in general, it is not possible to remove cubic terms and keep differentiable equivalence, though all higher order terms can be removed \cite{takens1973normal,glendinning2022normal}. 
We believe higher-order terms will be relevant for unfolding degenerate cases, such as where $\alpha$ or $\beta$ are zero. 

Our analysis of the local normal form bifurcations will have many implications on the global dynamics for a flow on a torus that contains such a singularity, including organising synchronisation and phase-locking properties of coupled oscillators near SNIC bifurcation (see, for example, \cite{ermentrout1986parabolic,laing2018dynamics}). In this sense, there is work still to do to generalise the results of \cite{baesens2013interaction} for the mixed case. There has recently been work \cite{augustsson2024coevolutionary} that considers global dynamics in such a case. However, they consider a specific family that does not have a parameter space of high enough dimension to unfold the local bifurcations in the sense we do here. In particular, that work considers a symmetric case which in our notation would have $\alpha=\beta$. Nonetheless, we note that the bifurcations they observe will be present in the more general case.

The current paper, as a study of quadratic planar vector fields, can be connected to Hilbert's 16th problem for the quadratic case. This concerns how many limit cycles are possible in a planar quadratic vector field, and remains unanswered \cite{dumortier1994hilberts,ilyashenko2002centennial}, even though there are only five independent parameters. The current maximum number of limit cycles found is four \cite{ilyashenko2002centennial,artes2021structurally}. In our study, we find regions with up to two limit cycles. However, we cannot rule out the possibility that there are additional saddle-node of limit cycle bifurcations on varying $\alpha$ and $\beta$ away from the chosen values. 

We do not address applications in this paper, but note that the problem that could be termed ``excitable competition", where multiple excitable connections exist from an attractor are needed to represent a node with multiple outgoing transitions in an excitable network \cite{ashwin2021excitable}. Such models are useful models of input-driven dynamics with finite-state computational properties; see, for example, \cite{ansmann2016selfinduced} or the review in \cite{ashwin2024}.

\subsection*{Acknowledgements}

 We thank LMS for support for CMP via a scheme 2 visiting lecturer grant. The work of CMP was also supported by the Royal
Society Te Aprangi (Marsden Fund Council, NZ Government), for funding via
Grant 21-UOA-048. For the purpose of open access, the authors have applied a Creative Commons Attribution (CC BY) licence to any Author Accepted Manuscript version arising from this submission.

 \subsection*{Data availability statement}

 The continuation of bifurcations was undertaken using COCO \cite{dankowicz2013recipes}. Additional code for plotting phase portraits and for continuing SNICeroclinic and non-central SNIC bifurcations (developed in Appendix~A) is available from \cite{SNinteraction-2025-pages}. 

\bibliographystyle{plainnat}

\bibliography{hilltoprefs}

\appendix

\newpage

\section{Defining systems for continuation of SNICeroclinic and Non-central SNIC bifurcations}
\label{sec:defining:systems}

The defining systems for non-central SNIC and SNICeroclinic, used to determine the curves in $(\mu,\beta,\gamma)$-space in Figure~\ref{fig:cod2-bifs}, follow the principles set out by \cite{beyn1990numerical}, and implemented in \textsc{HomCont} as part of \textsc{Auto} \cite{Doed07,Doed99}. The particular degeneracies corresponding to a non-central SNIC or a SNICeroclinic were not discussed in \cite{beyn1990numerical}. Hence, we specify them here. Our implementation makes a few modifications that are specific to the case of two-dimensional systems, and that simplify the defining system if one relies on finite differences of the Jacobian for the defining system. Our ODE system has the form written down in \eqref{eq:odesmu} for $u=(x,y)\in\R^2$, $p=(\mu,\beta,\gamma)\in\R^3$,
\begin{align}
    \label{gen:ode}
    \dot u&=f(u,p)\mbox{,}&\mbox{where\quad} 
    f((x,y),(\mu,\beta,\gamma))&=\begin{bmatrix}
        -\mu -\gamma+2\,\alpha \,y +x^2-\lambda_{{\Tr}0}\\
        \phantom{-}\gamma -\mu +2\,\beta \,x+y^2-\lambda_{{\Tr}0} 
    \end{bmatrix}
\end{align}
with $\lambda_{{\Tr}0}=\gamma(\gamma+\alpha^2-\beta^2)/(\alpha+\beta)^2$.
The systems and runs described below are implemented in \textsc{coco} \cite{dankowicz2013recipes}. Code 
is available at \url{https://github.com/jansieber/sniceroclinic-coco}.

\paragraph{ODE orbit segments}
Numerical boundary-value problem (BVP) solvers for ODEs, such as implemented in \textsc{coco}'s \textsc{coll} toolbox \cite{dankowicz2013recipes},
implement numerical approximations for general autonomous ODE orbit segments $u:[0,1]\to\R^{n_u}$ of the form
\begin{align}\label{gen:orbitseg}
    \dot u(t)&=Tf_\mathrm{ode}(u(t),p),& u(0)&=u_-,&u(1)&=u_+,    
\end{align}
where the variables $u_-,u_+\in\R^{n_u}$, $T\in\R$, $p\in\R^{n_p}$ can be fixed or varied as part of a continuation problem. 
Denoting the flow generated by the initial-value problem $\dot u=f_\mathrm{ode}(u,p)$, $u(0)=u_\mathrm{ini}$ as $\Phi_\mathrm{ode}(T;u_0,p)$, the orbit segment  equation \eqref{gen:orbitseg} is mathematically equivalent to the $n_u$ equations
    $0=\Phi_\mathrm{ode}(T;u_-,p)-u_+$.    
In our case we will have $n_u=2$ and $n_p=3$, since $u=(x,y)$ and $p=(\mu,\beta,\gamma)$.

\paragraph{Detection and continuation of non-central SNIC}

For a SNIC orbit we formulate the following system of $9$ equations for the $12$ variables $u_-\in\R^2$, $u_+\in\R^2$, $u_\mathrm{sn}\in\R^2$, $p=(\mu,\beta,\gamma)\in\R^3$, $s_-,s_+,T\in\R$:
\begin{align}\allowdisplaybreaks
\label{gen:ncsnic:orbit}
    \mbox{orbit ($2$ eqs.):}&&\dot u(t)&=Tf(u(t),p),\ u(0)=u_-, u(1)=u_+,\\
\label{gen:ncsnic:eq}
    \mbox{equilibrium ($2$ eqs.):}&&0&=f(u_\mathrm{sn},p),\\
\label{gen:ncsnic:sn}
    \mbox{saddle node ($1$ eq.):}&&0&=\det \partial_1f(u_\mathrm{sn},p),\\
\label{gen:ncsnic:bc-}
    \mbox{b.c. at $-\infty$ ($2$ eqs.):}&&0&=u_\mathrm{sn}-u_-+v_\mathrm{c}(u_\mathrm{sn},p,u_-)s_-,\\
\label{gen:ncsnic:bc+}
    \mbox{b.c. at $+\infty$ ($1$ eq.)\phantom{s}:}&&0&=w_\mathrm{c}(u_\mathrm{sn},p,u_-)^\tran[u_+-u_\mathrm{sn}]-s_+,\\
    \label{gen:ncsnic:phase}
    \mbox{phase cond. ($1$ eq.):}&&0&=\int_0^1u_\mathrm{r}(t)^\tran u(t)\mathrm{d} t+u_+^\tran\left[\frac{u_+}{2}-u_{\mathrm{r},+}\right]-u_-^\tran\left[\frac{u_-}{2}-u_{\mathrm{r},-}\right]. 
\end{align}
In \eqref{gen:ncsnic:phase} the variables $[0,1]\ni t\mapsto u_\mathrm{r}(t)\in\R^{n_u}$ and $u_{\mathrm{r},\pm}\in\R^{n_u}$ are reference solutions. At the initial step of the continuation they are equal to the the initial guess, and for all subsequent steps that are equal to the solution $[0,1]\ni t\mapsto u(t)\in\R^{n_u}$ and $u_\pm\in\R^{n_u}$ at a near-by previously computed point. System~\eqref{gen:ncsnic:orbit}--\eqref{gen:ncsnic:phase} assumes that the saddle node $u_\mathrm{sn}$ is transversally stable.

The functions $v_\mathrm{c}(u,p,u_\mathrm{dev})$ and  $w_\mathrm{c}(u,p,u_\mathrm{dev})$ in \eqref{gen:ncsnic:bc-} and \eqref{gen:ncsnic:bc+} are the right and left nullvectors $v$ and $w$ of $\partial_1f(u,p)$, made unique by scaling and adjusting orientation of $v$ as follows:
\begin{align}\label{gen:ncsnic:nullvecs}
    0&=\partial_1f(u,p)v,&
    0&=v^\tran v-1,&
    0&<v^\tran[u_\mathrm{dev}-u],&
    0&=\partial_1f(u,p)^\tran w,&
    0&=w^\tran v-1.
\end{align}

\begin{figure}
    \centering
    \includegraphics[width=0.5\linewidth]{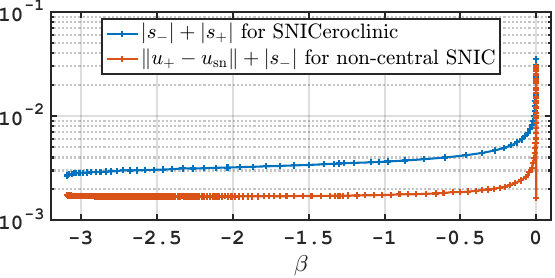}
    \caption{Distances between end of approximate numerical orbit segment and equilibrium along approximate SNICeroclinic and non-central SNIC shown in Figure~\ref{fig:cod2-bifs}, parametrized in $\beta$.}
    \label{fig:snicsdist}
\end{figure}
The initial values for the orbit segment $u(\cdot)$ and all variables are extracted from a periodic orbit of moderately large period that is near a SNIC. For a typical SNIC approximation $s_\pm$ will be small and satisfy $s_+<0<s_-$ by definition of $v_\mathrm{c}$, $w_\mathrm{c}$, and equations \eqref{gen:ncsnic:bc-} and \eqref{gen:ncsnic:bc+}. To obtain a more accurate approximation of a SNIC, to track the SNIC detecting a non-central SNIC, and finally to track the non-central SNIC, we performed $3$ runs.
\begin{description}
    \item[(Improve accuracy)] We vary $u_+$, $u_-$, $u_\mathrm{sn}$, $T$, $s_-$, $s_+$, $\mu$, and keep $\beta$, $\gamma$ fixed for dimension deficit equal to $1$, continuing a curve until we detect $T=2\times10^3$. During this run, parameters will not change; only the orbit segment will grow longer, increasing period $T$, and decreasing $|s_\pm|$.
    \item[(Track SNIC)]  We vary $u_+$, $u_-$, $u_\mathrm{sn}$, $s_-$, $s_+$, $\mu$, $\gamma$, and keep $\beta$, $T$ fixed for dimension deficit equal to $1$, continuing a curve until we detect $s_+=0$, indicating a non-central SNIC.
    \item[(Track non-central SNIC)]  We vary $u_+$, $u_-$, $u_\mathrm{sn}$, $s_-$, $\mu$, $\beta$, $\gamma$, and keep $s_+=0$ and $T$ fixed for dimension deficit equal to $1$, monitoring the distance 
    \begin{align*}
        d=\|u_--u_\mathrm{sn}\|+\|u_+-u_\mathrm{sn}\|=|s_-|+\|u_+-u_\mathrm{sn}\|,
    \end{align*}
    plotted in Figure~\ref{fig:snicsdist}. As long as this distance is small the result is a good approximation of a non-central SNIC. Note that the component $\|u_+-u_\mathrm{sn}\|$ is always less than $0.5\times 10^{-3}$, and that the deviation is of $u_-$ from $u_\mathrm{sn}$ is along the nullvector $v_\mathrm{c}$, which is tangent to the true trajectory (Figure~\ref{fig:snicsdist} does not display both parts separately).
\end{description}
\paragraph{Detection and continuation of SNICeroclinic}
For a SNICeroclinic orbit we formulate the following system of $12$ equations for the $14$ variables $u_-\in\R^2$, $u_+\in\R^2$, $u_\mathrm{sn}\in\R^2$, $u_\mathrm{sa}\in\R^2$, $p=(\mu,\beta,\gamma)\in\R^3$, $s_-,s_+,T\in\R$:
\begin{align}
\label{gen:homsnic:orbit}
    \mbox{orbit ($2$ eqs.):}&&\dot u(t)&=Tf(u(t),p),\ u(0)=u_-, u(1)=u_+,\\
\label{gen:homsnic:eq1}
    \mbox{eq. at $-\infty$ ($2$ eqs.):}&&0&=f(u_\mathrm{sn},p),\\
\label{gen:homsnic:sn}
    \mbox{saddle node ($1$ eq.):}&&0&=\det \partial_1f(u_\mathrm{sn},p),\\
\label{gen:homsnic:eq2}
    \mbox{eq. at $+\infty$ ($2$ eqs.):}&&0&=f(u_\mathrm{sa},p),\\
\label{gen:homsnic:bc-}
    \mbox{b.c. at $-\infty$ ($2$ eqs.):}&&0&=u_\mathrm{sn}-u_-+v_\mathrm{c}(u_\mathrm{sn},p,u_-)s_-,\\
\label{gen:homsnic:bc+}
    \mbox{b.c. at $+\infty$ ($2$ eqs.):}&&0&=u_\mathrm{sa}-u_++v_\mathrm{s}(u_\mathrm{sa},p,u_+)s_+,\\
    \label{gen:homsnic:phase}
    \mbox{phase cond. ($1$ eq.):}&&0&=\int_0^1u_\mathrm{r}(t)^\tran u(t)\mathrm{d} t+u_+^\tran\left[\frac{u_+}{2}-u_{\mathrm{r},+}\right]-u_-^\tran\left[\frac{u_-}{2}-u_{\mathrm{r},-}\right]. 
\end{align}
In \eqref{gen:homsnic:phase} the variables $ u_\mathrm{r}$ and $u_{\mathrm{r},\pm}$, and the function $v_\mathrm{c}(u,p,u_\mathrm{dev})$  have the same meaning as in \eqref{gen:ncsnic:bc-} and \eqref{gen:ncsnic:phase}, and are defined as in \eqref{gen:ncsnic:nullvecs}. The vector $v_\mathrm{s}(u,p,u_\mathrm{dev})$ is the eigenvector for the stable eigenvalue of the saddle point $u_\mathrm{sa}$, made unique by scaling and adjusting orientation. It equals the vector $v$ satisfying the conditions
\begin{align}
    \label{gen:homsnic:evec}
    0&=\partial_1f(u,p)v-\lambda v,&0&=v^\tran v-1,&\lambda&<0,&0&<v^\tran[u_\mathrm{dev}-u].
\end{align}
System~\eqref{gen:homsnic:orbit}--\eqref{gen:homsnic:phase} assumes that the saddle node $u_\mathrm{sn}$ is transversally stable. The initial values for the orbit segment $u(\cdot)$ and all variables are extracted from a periodic orbit of moderately large period that is near a SNICeroclinic. Along the branch of long-period periodic orbits shown in blue in Figure~\ref{fig:bifset-abn}, a large number of steps is spent near point (2) (black triangle pointing upwards), which is the approximate SNICeroclinic. These periodic orbits spend two long subperiods near equilibria, one of them is approximately $u_\mathrm{sn}$, the other is $u_\mathrm{sa}$. For a typical SNICeroclinic approximation $s_\pm$ will be small and satisfy $0<s_+$ and $0<s_-$ by definition of $v_\mathrm{c}$, $v_\mathrm{s}$, and equations \eqref{gen:homsnic:bc-} and \eqref{gen:homsnic:bc+}. To obtain a more accurate approximation of the SNICeroclinic, and then track the SNICeroclinic we perform $2$ runs.
\begin{description}
    \item[(Improve accuracy)] We vary $u_+$, $u_-$, $u_\mathrm{sn}$, $u_\mathrm{sa}$, $T$, $s_-$, $s_+$, $\mu$, $\gamma$ and keep $\beta$ fixed for dimension deficit equal to $1$, continuing a curve until we detect $T=10^3$. During this run, parameters do not change but the orbit segment grows longer, increasing period $T$, and decreasing $s_\pm$.
    \item[(Track SNICeroclinic)]  We vary $u_+$, $u_-$, $u_\mathrm{sn}$, $u_\mathrm{sa}$, $s_-$, $s_+$, $\mu$, $\beta$, $\gamma$, and keep $T$ fixed for dimension deficit equal to $1$, monitoring the distance $d=s_-+s_+$,
    plotted in Figure~\ref{fig:snicsdist}. As long as this distance is small the result is a good approximation of a SNICeroclinic. Note that the component $s_2$ is always less than $0.52\times 10^{-4}$, such that the approximation accuracy is limited mostly by the gap between $u_\mathrm{sn}$ and $u_-=u(0)$.
\end{description}

\end{document}